\documentclass{amsart}

\usepackage{url, hyperref}
\usepackage{cite}

\usepackage{amssymb}
\usepackage{bm}
\usepackage{graphicx}
\usepackage{psfrag}
\usepackage{color}
\usepackage{url}
\usepackage{amsmath,color}
\usepackage{enumerate}
\usepackage{algpseudocode} 
\usepackage{lmodern}
\usepackage{todonotes}
\usepackage{xy}
\input xy
\xyoption{all}
\usepackage{mathrsfs}

\newcommand{\excise}[1]{}

\newcommand{\ve}{\boldsymbol}

\newtheorem{thm}{Theorem}
\newtheorem{lemma}[thm]{Lemma}
\newtheorem{cor}[thm]{Corollary}
\newtheorem{prop}[thm]{Proposition}

\newtheorem{remark}{Remark}

\newcommand{\be}{\begin{eqnarray}}
\newcommand{\bea}{\begin{eqnarray*}}
\newcommand{\ee}{\end{eqnarray}}
\newcommand{\eea}{\end{eqnarray*}}

\newcommand{\ring}[1]{\ensuremath{\mathbb{#1}}}

\newcommand{\diag}{\mathrm{diag}}

\newcommand{\lin}{\mathrm{lin}}

\newcommand{\conv}{\mathrm{conv}}
\newcommand{\vol}{\mathrm{vol}}

\newcommand{\R}{\mathbb R}

\newcommand{\Z}{\mathbb Z}

\newcommand{\KP}{{P}}

\newcommand{\dl}{l}
\newcommand{\ub}{f}

\newcommand{\PF}{{p}}

\newcommand\CP{{C}}

\newcommand\Q{\ring{Q}}

\newcommand\supp{\mathrm{supp}}

\newcommand\bb{\ve b}%
\newcommand\cc{\ve c}%
\newcommand\xx{\ve x}%

\newcommand\BB{\ve B}%
\renewcommand\AA{\ve A}%
\renewcommand\SS{\ve S}
\newcommand\ID{\ve I}

\newcommand\DD{\ve D}%

\newcommand\UU{\ve U}%

\def\udots{\mathinner{\mkern1mu\raise1pt\vbox{\kern7pt\hbox{.}}\mkern2mu
\raise4pt\hbox{.}\mkern2mu\raise7pt\hbox{.}\mkern1mu}}

\global\long\def\zero{\mathbf{0}}

\title[Optimal Proximity Bound and Product Function Estimates]{Optimal Proximity Bound and Product Function Estimates in Integer Linear Programming}

\author{Iskander Aliev}

\address{Cardiff University, alievi@cardiff.ac.uk}

\author{Gennadiy Averkov}

\address{Brandenburg University of Technology Cottbus-Senftenberg, averkov@b-tu.de}

\author{William Jones}

\address{Cardiff University, jonesw19@cardiff.ac.uk}

\author{Timm Oertel}

\address{Friedrich-Alexander-Universität Erlangen–Nürnberg, timm.oertel@fau.de}

\date{\today }

\begin{document}
\maketitle
\thispagestyle{plain}

\begin{abstract}
	\noindent
We obtain an optimal proximity bound for integer linear programs in standard form
\begin{eqnarray}\label{abstract}
\max\{\cc\cdot\xx:\AA\xx=\bb, \xx\in \Z^{n}_{\ge 0}\}\,,
\end{eqnarray}
where $\AA$ is an integer $m\times n$ matrix of rank $m<n$ and $\bb$ is an integer vector. Specifically, we show that the Euclidean distance from any optimal vertex solution of the LP relaxation of \eqref{abstract} to a nearest optimal integer solution is bounded by $\sqrt{\det(\AA\AA^\top)}-1$ and that this estimate is asymptotically tight. 
We also derive bounds for the optimal solutions to \eqref{abstract} involving the product function $\prod_{i=1}^{n}(x_i+1)$ and discuss their applications in the knapsack setting.

\end{abstract}

\section{Introduction and main results}

Let $\AA \in \Z^{m\times n}$ with $m\le n$ and $\ve b \in \Z^m$.
We define the polyhedron
\bea
\KP(\AA, {\ve b})=\{{\ve x}\in \R^n_{\ge 0}: \AA{\ve x}={\ve b}\}\,
\eea
and assume that $\KP(\AA, {\ve b})$ contains integer points. 
Given a set $K\subset \R^n$, we denote by $\conv(K)$ its convex hull. The polyhedron 
\bea
\KP_I(\AA,\bb)=\conv(\KP(\AA,\bb)\cap \Z^n)\,
\eea
is traditionally referred to as the {\em integer hull} of $\KP(\AA,\bb)$. 

We consider the integer linear programming
problem 
\begin{eqnarray}
\max\{\cc\cdot\xx:\xx\in\KP(\AA,\bb)\cap\Z^{n}\}\,,\label{initial_IP}
\end{eqnarray}
where $\cc\in\Q^{n}$ is a cost vector and the $\cc\cdot\xx$ denotes the standard inner product. Throughout the paper, we assume that (\ref{initial_IP}) is feasible and bounded. Note that at least one  optimal solution of \eqref{initial_IP} is a vertex of the integer hull $\KP_I(\AA,\bb)$.

One of the classical approaches for solving optimisation problems 
of the form \eqref{initial_IP} is based on solving the
linear programming relaxation 
\begin{eqnarray}
\max\{\cc\cdot\xx:\xx\in\KP(\AA,\bb)\}\,,\label{initial_LP}
\end{eqnarray}
obtained
by dropping the integrality constraint. Various techniques are then used to construct an optimal integer solution ${\ve z}^*$  of \eqref{initial_IP} from a fractional solution ${\ve x}^*$ of \eqref{initial_LP}. In this setting, proximity-type results play a central role by providing bounds on the distance between ${\ve x}^*$ and ${\ve z}^*$ in terms of structural parameters of the matrix $\AA$. 

We now introduce notation needed to state the relevant proximity results. For $\tau =\{i_1, \ldots, i_k\}\in \binom{[n]}{k}$, let $\AA_\tau$ denote the $m \times k$ submatrix of $\AA$ with columns indexed by $\tau$. Similarly, for ${\ve x}\in \R^n$, we denote by ${\ve x}_\tau$ the subvector $(x_{i_1}, \ldots, x_{i_k})^{\top}$. For $1\le r\le m$,  we denote by $\Delta_r(\AA)$ the maximum absolute $r\times r$ subdeterminant of $A$, that is 
\bea
\Delta_r(\AA)=\max\left\{|\det \BB|: \BB \mbox{ is an }r\times r \mbox{ submatrix of }\AA\right\}\,.
\eea
Further, $\|\cdot\|_\infty$,  $\|\cdot\|_1$ and $\|\cdot\|_2$ denote the $\ell_\infty$, $\ell_1$ and $\ell_2$ norms, respectively.

Let ${\ve x}^*$ be a vertex optimal solution of \eqref{initial_LP}.
The sensitivity theorems of Cook et al. 
\cite{MR839604}  imply the existence of an optimal
solution ${\ve z}^*$ to \eqref{initial_IP} satisfying
\be\label{PROX_C}
\|{\ve x}^*-{\ve z}^* \|_\infty\le (n-m) \Delta_m(\AA)\,.
\ee
A recent improvement  by Celaya et al. \cite{CKPW,CKPW1} refines  \eqref{PROX_C} to
\be\label{PROX_CKPW}
\|{\ve x}^*-{\ve z}^* \|_\infty\le \frac{4(n-m)+2}{9} \Delta_m(\AA)\,.
\ee
 In the same setting,  Eisenbrand and Weismantel \cite{MR3775840} obtain the bound
\be\label{PROX_EW}
\|{\ve x}^*-{\ve z}^* \|_1\le m (2m\Delta_1(\AA)+1)^m\,.
\ee
For other contributions, we refer the reader to Celaya and Henk \cite{CH}, Lee et al. \cite{LPSX,LPSX1},  Oertel, Paat and Weismantel \cite{OPW}, Paat, Weismantel and Weltge \cite{PWW}, as well as the references therein.

We will  use the quantity
\bea
\Delta(\AA)=\sqrt{\det \AA\AA^\top}\,.
\eea
Geometrically, $\Delta(\AA)$ is the $m$-dimensional volume of the parallelepiped determined by the rows of $\AA$.
The Cauchy-Binet formula implies
\be\label{CB}
\begin{aligned}
\Delta(\AA) &=\sqrt{\sum_{\substack{\sigma\subset[n]: \\ |\sigma|=m}}(\det(\AA_\sigma))^2}\,.
\end{aligned}
\ee
We also denote by $\gcd(\AA)$ the greatest common divisor of all $m\times m$ subdeterminants of $\AA$. 

The first result of this paper establishes an optimal proximity bound for the Euclidean distance between an optimal vertex solution
${\ve x}^*$ of the linear programming relaxation \eqref{initial_LP}  and its nearest optimal solution ${\ve z}^*$ of the original integer program \eqref{initial_IP} in terms of the quantity $\Delta(\AA)$.
To the best of our knowledge, this is the first optimal general proximity bound known in the literature.

%
\begin{thm}
\label{thm:UpperBound}
Let $\AA \in \Z^{m\times n}$, $m<n$, be a matrix of rank $m$, $\ve b \in \Z^m $, and suppose that $ \KP(\AA, {\ve b})$ contains integer points.  Let ${\ve x}^*$ be an optimal vertex solution to \eqref{initial_LP}. Then there exists a vertex ${\ve z}^*$ of $\KP_I(\AA, {\ve b})$, such that ${\ve z}^*$  is an optimal solution to \eqref{initial_IP}  and the bound
\begin{equation}\label{t:Bound_via_pencils}
\|{\ve x}^*-{\ve z}^*\|_2\le  \frac{\Delta(\AA)}{\gcd(\AA)}-1\,
\end{equation}
holds.
\end{thm}
We remark that the bound \eqref{t:Bound_via_pencils} was known in the feasibility setting; see Theorem 3 in \cite{ACH_IPCO,ACH_MPB}.
For certain matrices $\AA$ the bound \eqref{t:Bound_via_pencils} is smaller than the $\ell_2$-norm proximity bounds $\sqrt{n}((4(n-m)+2)/9)\Delta_m(\AA)$  and $m (2m\Delta_1(\AA)+1)^m$ that follow from \eqref{PROX_CKPW} and \eqref{PROX_EW}, respectively.  Indeed, it suffices to observe that the ratio $\Delta(\AA)/\Delta_m(\AA)$ can be arbitrarily close to one and the ratio  $\Delta(\AA)/\Delta_1(\AA)$ can be arbitrarily small. For instance, consider the sequence of matrices $\AA_t=(t{\ve I}_m|\BB)$, where ${\ve I}_m$ is the $m\times m$ identity matrix and  $\BB$ is a fixed $m\times(n-m)$ integer matrix. Then  $\Delta(\AA_t)/\Delta_m(\AA_t)\rightarrow 1$ as $t\rightarrow\infty$. Furthermore, the quantity $\Delta(\AA)$ equals the determinant of the lattice generated by the rows of $\AA$. The same lattice can be generated by the rows of a matrix $\widetilde\AA$ with arbitrarily large $\Delta_1(\widetilde\AA)$. Hence, the ratio $\Delta(\widetilde\AA)/\Delta_1(\widetilde\AA)$ can be made arbitrarily small. 

Next, we show that the proximity bound \eqref{t:Bound_via_pencils} is optimal.

\begin{thm}
\label{thm:UpperBoundOptimality}
For any integers  $1\le m<n$ and for any $\epsilon>0$, there exist a matrix $\AA \in \Z^{m\times n}$ of rank $m$, $\ve b \in \Z^m $ and ${\ve c}\in \Z^n$ with the following property. The linear relaxation \eqref{initial_LP} and the original problem \eqref{initial_IP} have unique optimal solutions ${\ve x}^*$ and ${\ve z}^*$, respectively,  with
\begin{equation}\label{t:Bound_via_pencils_optimality}
\|{\ve x}^*-{\ve z}^*\|_2> (1-\epsilon)  \left(\frac{\Delta(\AA)}{\gcd(\AA)}-1\right)\,.
\end{equation}
\end{thm}

We now consider the product  function $\PF({\ve x}): \R^{n}\rightarrow \R_{\ge 1}$ defined as 
\bea
\PF({\ve x})=\prod_{i=1}^{n}(x_i+1)\,.
\eea

The use of the  function $\PF({\ve x})$ in integer programming dates back to a seminal paper by Gomory \cite{Gomory_polyhedra}, where it was applied in the context of asymptotic integer programming and corner polyhedra. In \cite{ACH_IPCO,ACH_MPB,ACHW}, this function was used to derive proximity-sparsity transference bounds for integer programs. Klein \cite{Klein} used the product function to obtain new algorithmic results for the {\em Unbounded SubsetSum}, {\em Frobenius}, and {\em All-Targets} problems. Deng, Mao and Zhong \cite{DMZ} used $\PF({\ve x})$ in their work on the {\em All-Target Unbounded Knapsack}, {\em All-Targets CoinChange}, and {\em All-Targets} problems.

We will first consider the problem of estimating the minimum value of  $\PF({\ve x})$ over optimal solutions to \eqref{initial_IP}. This question was first studied in \cite{DMZ} in the unbounded knapsack scenario $\AA\in\Z^{1\times n}_{> 0}$. Specifically, Lemma 3.3 in \cite{DMZ} implies in this setting that the problem \eqref{initial_IP} has an optimal solution ${\ve z}^*$ with 
\be\label{DMZ-bound}
\PF({\ve z}^*)\le b+1= \AA{\ve z}^*+1\,.
\ee
The bound \eqref{DMZ-bound} is tight (see Remark \ref{tightness_b_plus_1} in the proof of Theorem \ref{thm_EL_proximity_strong_ordered_optimality} below).

In this paper, we obtain the first upper bound on the product function that applies in the general case. For ${\ve x}\in \R^n$, we let $\DD({\ve x})=\diag(x_1+1, \ldots, x_n+1)$ and 
\bea
\ub_{\AA}({\ve x})=\frac{\Delta(\AA\DD({\ve x}))}{\gcd(\AA)}\,.
\eea
\begin{thm}\label{thm_EL_interger_hull}
Let $\AA\in \Z^{m\times n}$, $m<n$,  be a matrix of rank $m$ and let  ${\ve b}\in \Z^m$. Assume that the polyhedron $\KP(\AA, {\ve b})$ contains integer points.  Then any  vertex ${\ve z}^*$ of $\KP_I(\AA, {\ve b})$ satisfies the bound
\be\label{EL_integer_hull_proximity_bound} 
\PF({\ve z}^*)\le \ub_{\AA}({\ve z}^*)\,.
\ee
\end{thm}

The proof of Theorem \ref{thm_EL_interger_hull} relies on results from Minkowski's geometry of numbers. This geometric approach naturally leads to the choice of the function $\ub_{\AA}({\ve z}^*)$ in \eqref{EL_integer_hull_proximity_bound}. Notice that, in the general case, a bound on the product function cannot depend only on ${\ve b}$ and dimensional constants. Already for unbounded knapsack problems with negative entries, $\PF({\ve z}^*)$ can be arbitrarily larger than $b$ (see Remark~\ref{worse_than_b} in the proof of Theorem~\ref{thm_EL_integer_hull_proximity_bound_optimality}). 

Using \eqref{CB}, we can write
\be\label{ub_expanded}
\begin{aligned}
\Delta(\AA\DD({\ve x})) &=\sqrt{\sum_{\substack{\sigma\subset[n]: \\ |\sigma|=m}}(\det(\AA_\sigma) \prod_{i\in \sigma}(x_i+1) )^2}\,.
\end{aligned}
\ee
Hence, $(\ub_{\AA}({\ve x}))^2$ is a polynomial of total degree $2m$ in ${\ve x}$, whereas the right hand side of \eqref{DMZ-bound} is the value of the total degree one polynomial $\AA{\ve x} +1$
at ${\ve x} = {\ve z}^*$. This suggests that \eqref{EL_integer_hull_proximity_bound} may be viewed as a natural generalisation of \eqref{DMZ-bound}.  

Observe that it is not difficult to construct examples where \eqref{EL_integer_hull_proximity_bound} is tight and ${\ve z}^*$ is a vertex of $\KP(\AA, {\ve b})$. For instance, let $\AA=(\ID_m | {\ve 0}_{m\times (n-m)})$, where $\ID_m$ denotes the $m\times m$
identity matrix and ${\ve 0}_{m\times (n-m)}$ denotes the zero $m\times (n-m)$ matrix. Let further ${\ve b}=(b_1, \ldots, b_m)^\top$ be any nonnegative integer vector. Then ${\ve z}^*=(b_1, \ldots, b_m, 0,\ldots, 0)^\top$ is a vertex of both $\KP(\AA, {\ve b})$ and $\KP_I(\AA, {\ve b})$, and 
\bea
\PF({\ve z}^*)= \ub_{\AA}({\ve z}^*)=\prod_{i=1}^m (b_i+1)\,.
\eea

In the next result, we exclude vertices of $\KP(\AA, {\ve b})$ and show that  the bound \eqref{EL_integer_hull_proximity_bound} remains optimal up to the multiplicative constant $2/\sqrt{5}$.
\begin{thm}
\label{thm_EL_integer_hull_proximity_bound_optimality}
For any integers $1\le m<n$ and for any $\epsilon>0$, there exist a matrix $\AA \in \Z^{m\times n}$ of rank $m$ and $\ve b \in \Z^m $ with the following property.  The integer hull $\KP_I(\AA, {\ve b})$ has a vertex ${\ve z}^*$ that is not a vertex of $\KP(\AA, {\ve b})$ such that
\begin{equation}\label{EL_integer_hull_proximity_optimality}
\PF({\ve z}^*)> \left (\frac{2}{\sqrt{5}}-\epsilon\right)  \ub_{\AA}({\ve z}^*)\,.
\end{equation}
\end{thm}

We now consider nonnegative matrices $\AA\in \Z^{m\times n}_{\ge 0}$. In this case, Theorem~\ref{thm_EL_interger_hull} implies a bound that depends only on ${\ve b}$ and $m$, similarly to~\eqref{DMZ-bound}. Given a vector ${\ve x}=(x_1,\ldots, x_n)^\top\in \R^n$, let $\supp({\ve x})$ denote the {\em support} of ${\ve x}$, that is $\supp({\ve x})=\{i: x_i\neq 0\}$. 
\begin{thm}\label{thm_nonnegative}
Let $\AA\in \Z^{m\times n}_{\ge 0}$, $m<n$,  be a matrix of rank $m$ and let  ${\ve b}\in \Z^m$. Assume that the polyhedron $\KP(\AA, {\ve b})$ contains integer points.  Then any  vertex ${\ve z}^*$ of $\KP_I(\AA, {\ve b})$ 
satisfies the bound
\be\label{nonnegative_A} 
\PF({\ve z}^*)\le 2^m \prod_{i\in\supp({\ve b})} b_i\,.
\ee
\end{thm}
Given ${\ve x}=(x_1,\ldots, x_n)^\top\in\Z^n_{\ge 0}$, we denote by  $\langle {\ve x} \rangle$ its {\em binary encoding length} defined as
\bea
\langle {\ve x} \rangle= \sum_{i=1}^n\lceil\log_2(x_i+1)\rceil\,.
\eea
Theorem \ref{thm_nonnegative} immediately implies a neat bound for the function $\langle {\ve x} \rangle$ in the case of nonnegative matrices $\AA$.
\begin{cor}\label{cor_nonnegative}
In the setting of Theorem \ref{thm_nonnegative}
\be\label{nonnegative_A} 
\langle {\ve z}^* \rangle < m+n+ \sum_{i\in\supp({\ve b})}  \log_2(b_i)\,. 
\ee
\end{cor}

Recall that $\gamma\subset[n]$ is a {\em basis} of $\AA$ if $|\gamma|=m$ and the submatrix $\AA_{\gamma}$  is nonsingular. For any vertex $\ve x^*$  of $\KP(\AA, {\ve b})$  %
there is a basis $\gamma$ of $A$ such that
\be\label{vertex_assumption}
 {\ve x}^*_\gamma = A_\gamma^{-1}{\ve b} \mbox{ and }{\ve x}^*_{\bar \gamma}={\ve 0}\,,
\ee
where ${\bar \gamma}$ denotes the complement of $\gamma$ in $[n]$.
In general, for a given vertex there can be many choices for the basis $\gamma$ in \eqref{vertex_assumption}.
The variables indexed by $\tau$ and $\bar\tau$ are often  referred to as {\em basic} and {\em nonbasic}  variables, respectively.

The traditional proximity results \eqref{PROX_C}-\eqref{t:Bound_via_pencils} give estimates for various norms of the difference vector ${\ve r}={\ve x}^*-{\ve z}^*$ in terms of the parameters of the matrix $\AA$. 
Observe that \eqref{vertex_assumption} implies 
\bea {\ve r}_\gamma=\AA_\gamma^{-1}{\ve b} -{\ve z}^*_\gamma \mbox{ and }{\ve r}_{\bar\gamma}=-{\ve z}^*_{\bar\gamma}\,.
\eea 
Further ${\ve z}^*_\gamma=\AA_\gamma^{-1} {\ve b}-\AA_\gamma^{-1}(\AA_{\bar\gamma} {\ve z}^*_{\bar\gamma})$.
Consequently, the difference vector ${\ve r}$ can be uniquely reconstructed from its projection ${\ve r}_{\bar\gamma}=-{\ve z}^*_{\bar\gamma}$ using the formula
\bea
{\ve r}= \left(
\begin{array}{c}
\AA_\gamma^{-1} (\AA_{\bar\gamma}{\ve z}^*_{\bar\gamma})\\
-{\ve z}^*_{\bar\gamma}
\end{array}
\right)\,.
\eea

Hence, from a computational standpoint, it is natural to consider projecting the solutions of \eqref{initial_IP} onto the nonbasic variables, with the goal of estimating the size of a region containing the subvector ${\ve r}_{\bar\gamma}=-{\ve z}^*_{\bar\gamma}$. In this setting, the choice of $\PF({\ve x})$ as a proximity function  already appears naturally  in the work of Gomory  \cite{Gomory_polyhedra} on group relaxations and corner polyhedra. 

Let $\ve x^*$  be a vertex of $\KP(\AA, {\ve b})$  with basis $\gamma$.  Following Gomory  \cite{Gomory_polyhedra}, we define  the {\em corner polyhedron} $\CP_{\gamma}(\AA, {\ve b})$ associated with $\gamma$ as
\bea
\CP_{\gamma}(\AA,\bb)=\conv(\{\xx\in\Z^{n}:\AA\xx=\bb,\xx_{\bar{\gamma}}\geq\zero\}).
\eea
Theorem 1 in \cite{Gomory_polyhedra} implies (see Theorem 11 in \cite{ACHW}) that 
any vertex ${\ve z}^{*}$ of the corner polyhedron $\CP_{\gamma}(\AA, {\ve b})$ 
satisfies the bound
\be\label{Gomory_corner}
\PF({\ve z}^*_{\bar\gamma})\leq \frac{|\det(\AA_{\gamma})|}{\gcd(\AA)}\,.
\ee
The bound \eqref{Gomory_corner} is particularly useful in the knapsack setting $\AA\in \Z^{1\times n}_{> 0}$, where it leads to the proximity-sparsity transference bounds from \cite{ACHW}. The main reason is that  at least one vertex of the integer hull $\KP_I(\AA, {\ve b})$ is also a vertex of
the corner polyhedron $\CP_{\gamma}(\AA,\bb)$. This property does not hold in general and already breaks for knapsacks with negative entries (see, e.g., Example 3.3.1 in \cite{Aled_thesis}).
This naturally leads to the question of estimating $\PF({\ve z}^*_{\bar\gamma})$ over the vertices ${\ve z}^*$ of $\KP_I(\AA, {\ve b})$
in the general case.

We first show that the vertices of the integer hull $\KP_I(\AA,\bb)$ possess ``small'' $(n-m)$-dimensional coordinate projections in terms of
the function $\PF({\ve x})$.

\begin{thm}\label{thm_EL_proximity_strong_ordered}
Let $\AA\in \Z^{m\times n}$, $m<n$,  be a matrix of rank $m$, and let  ${\ve b}\in \Z^m$.  Assume that $\KP(\AA, {\ve b})$ contains integer points and let  ${\ve z}^*$ be  a vertex of $\KP_I(\AA, {\ve b})$ with $z^*_{i_1}\le z^*_{i_2}\le\cdots\le z^*_{i_n}$.  Letting $\tau=\{i_{1}, \ldots, i_{n-m}\}$, we have %
\be\label{HW_proximity_bound_order} 
\PF({\ve z}^*_{\tau})\le  \frac{\Delta(\AA)}{\gcd(\AA)} \,.
\ee
\end{thm}

Next we show that the bound \eqref{HW_proximity_bound_order} is optimal up to the multiplicative constant $\sqrt{3}/{2}$. 
\begin{thm}
\label{thm_EL_proximity_strong_ordered_optimality}
For any positive integer $m$ and for any $\epsilon>0$, there exist a matrix $\AA \in \Z^{m\times n}$ of rank $m<n$ and $\ve b \in \Z^m $with the following property.  The integer hull $\KP_I(\AA, {\ve b})$ has a vertex ${\ve z}^*$ such that, using the notation from Theorem \ref{thm_EL_proximity_strong_ordered}, 
\begin{equation}\label{HW_proximity_bound_order_optimality}
\PF({\ve z}^*_{\tau})> \left (\frac{\sqrt{3}}{2}-\epsilon\right)  \frac{\Delta(\AA)}{\gcd(\AA)}\,.
\end{equation}
\end{thm}

We derive two corollaries of Theorem \ref{thm_EL_proximity_strong_ordered}.
To measure the size of the support of a vector ${\ve x}\in\R^n$, we use the {\em $0$-norm} $\|{\ve x}\|_0=|\supp({\ve x})|$, which is widely used in {compressed sensing} \cite{candestao} and integer optimisation \cite{AADO,ACH_IPCO,ACH_MPB,ACHW, Support, ADON2017}.

Theorem 1 in \cite{Support} gives an upper bound on the support size of optimal solutions to general integer programs. This result, together with its proof, implies that  every vertex of $\ve z^*$ of the integer hull $\KP_I(\AA, {\ve b})$ satisfies 
\begin{equation}\label{t:Bound_via_determinant_bound1}
\|\ve z^*\|_0 \le m + \log_2\left(\frac{\Delta(\AA)}{\gcd(\AA)}\right)\,.
\end{equation}
Theorem 2 in \cite{Support} and Lemma 15 in \cite{BJK} show that the bound \eqref{t:Bound_via_determinant_bound1} is close to optimal.
However, the question of its optimality remains open. The first corollary of Theorem \ref{thm_EL_proximity_strong_ordered} provides a refined sparsity bound that additionally takes into account the $n-m$ smallest entries of  ${\ve z^*}$. 
\begin{cor}\label{cor_general_sparsity}
Under the assumptions of Theorem  \ref{thm_EL_proximity_strong_ordered}, we have
\be\label{delicate_general_sparsity} 
\|{\ve z}^*\|_0\le  m+ \log_2\left(\frac{2^{\|{\ve z}^*_{\tau}\|_0}\Delta(\AA)}{\PF({\ve z}^*_\tau)\gcd(\AA)} \right)\,.
\ee
\end{cor}
The bound \eqref{delicate_general_sparsity} strengthens \eqref{t:Bound_via_determinant_bound1} unless 
${\ve z}^*_{\tau}\in \{0,1\}^{n-m}$.

One useful property of the function $\PF({\ve x})$ is that the ratio
$\|{\ve x}\|_1/\PF({\ve x})$ is exponentially small as a function of $\|{\ve x}\|_0$. This yields the following corollary.

\begin{cor}\label{cor_EL_proximity_strong_ordered}
Suppose that the vertex  ${\ve z}^*$ of $\KP_I(\AA, {\ve b})$ in the statement of Theorem \ref{thm_EL_proximity_strong_ordered} has $d=\|{\ve z}^*_{\tau}\|_0\ge 2$. Then the bound 
\be\label{cor_proximity_bound_order} 
\|{\ve z}^*_{\tau}\|_1\le \frac{d}{2^d}\left\lfloor\frac{\Delta(\AA)}{\gcd(\AA)}\right\rfloor
\ee
holds.
\end{cor}
We discuss algorithmic implications of Corollary \ref{cor_EL_proximity_strong_ordered} for  knapsack problems with sufficiently {\it dense} solutions in Section \ref{algorithmic implications}.

Next, we combine Theorems \ref{thm:UpperBound} and \ref{thm_EL_proximity_strong_ordered} to estimate the size of a region containing the subvector ${\ve r}_{\bar\gamma}=-{\ve z}^*_{\bar\gamma}$ in terms of  $\PF({\ve z}^*_{\bar\gamma})$ and $\Delta(\AA)$. For simplicity, we assume without loss of generality that $\gcd(\AA)=1$.

\begin{thm}\label{thm_EL_proximity_strong}
Let $\AA\in \Z^{m\times n}$, $m<n$,  be a matrix of rank $m$ with $\gcd(\AA)=1$, let  ${\ve b}\in \Z^m$ and  ${\ve c}\in \Z^n$. Assume that the problem \eqref{initial_IP} is feasible and bounded. Let $\ve x^*$ be a  vertex optimal solution to \eqref{initial_LP} with a basis $\gamma$. Then there exists a vertex ${\ve z}^*$ of $\KP_I(\AA, {\ve b})$, such that ${\ve z}^*$  is an optimal solution to \eqref{initial_IP}  and 
\bea 
\PF({\ve z}^*_{\bar\gamma})\le \Delta(\AA)\left(\frac{\Delta(\AA)-1}{\sqrt{m}}+1\right)^m\,.
\eea
\end{thm}

Compared with the bound \eqref{HW_proximity_bound_order}, fixing a basis $\gamma$ in Theorem \ref{thm_EL_proximity_strong} results in an additional factor of order $\Delta(\AA)^m$. Nevertheless, when $m$ is a fixed small integer and the dimension $n$ is large, the bound looks reasonably strong. It is unlikely, however, that the bound is close to optimal.

\section{Proximity function $\PF(\cdot)$ and knapsack problems with dense solutions}
\label{algorithmic implications}

The algorithmic implications of general proximity bounds involving the product function $\PF({\ve x})$ are a topic for future research. In this section, we focus on the unbounded knapsack scenario $\AA \in \Z^{1\times n}_{>0}$. For convenience, we avoid  matrix notation in this case and write $\AA={\ve a}^\top$ with ${\ve a}\in\Z^n_{>0}$. 
We will consider the {\em knapsack polytope}
\bea
\KP({\ve a}, {b})=\{{\ve x}\in \R^n_{\ge 0}: {\ve a}\cdot {\ve x}={b}\}\,,
\eea
where ${b}$ is a positive integer. Then the problem \eqref{initial_IP}, referred to as the {Unbounded Knapsack Problem}, takes the form
\be \label{initial_IP_knapsack} 
\max\{\cc\cdot\xx:\xx\in\KP({\ve a},b)\cap\Z^{n}\}\,.
\ee

A classical dynamic programming approach  for solving the knapsack problem \eqref{initial_IP_knapsack}  is based on constructing a directed weighted graph $G=(V, A)$ with vertex set $V\subset \Z_{\ge 0}$ and  representing  optimal solutions to \eqref{initial_IP_knapsack} as maximum-cost paths from $0$ to $b$.  For further details, we refer the reader to Schrijver \cite[Section 18.5]{Schrijver}, Papadimitriou \cite{Papadimitriou}, Eisenbrand and Weismantel \cite{MR3775840} and the references within.

Denote by  $D[t]$ the maximum in the problem \eqref{initial_IP_knapsack} with right-hand side $b=t$, and by $D[t_1,\ldots, t_k]$ the array consisting of the values $D[t_i]$. For convenience, we  also write $u=\|{\ve a}\|_\infty$. Interestingly, several dynamic programming algorithms for solving \eqref{initial_IP_knapsack}, or its inequality form, in fact solve the problem for a range of the values of $b$. For instance, the classical algorithm of Gilmore and Gomory (see, e. g., \cite{Gilmore_Gomory} and \cite[Section 3.6.2]{Martello_Toth}) solves \eqref{initial_IP_knapsack} for all targets $\le b$. Another example is a more recent algorithm due to Chan and He \cite{ChanHe}. Although their algorithm is designed for  the inequality form of the knapsack problem,  it can be  modified in a straightforward way to solve the equality form \eqref{initial_IP_knapsack}. The algorithm computes a collection of arrays $D[t-u, \ldots, t]$.

In this section, we assume that we are given a dynamic programming algorithm that computes the array $D[t-u+1, \ldots, t]$. Equivalently, we assume that  solutions to \eqref{initial_IP_knapsack} are available for $u$ consecutive right-hand sides $t-u+1, \ldots, t$.

Note that the size of the graph $G$ grows with $b$, which may be arbitrarily large. To tackle this issue, the dynamic programming algorithms in \cite{ChanHe, DMZ} use a preprocessing step that reduces $b$ below $u^2$ by subtracting a suitable multiple of an item with maximum value-to-weight ratio.

Perhaps surprisingly, Corollary \ref{cor_EL_proximity_strong_ordered} shows that $b$, and consequently the size of $G$, can be reduced significantly under the assumption that there exists an optimal vertex ${\ve z}^*$ of the integer hull $\KP_I({\ve a},b)$ with relatively large size of the support. We will refer to such optimal solutions as {\it dense}. More precisely, we will show that, assuming that $\|{\ve z}^*\|_0>\dl$ for a sufficiently large density level $\dl$, the right-hand side ${b}$ can be reduced to a quantity exponentially smaller in $\dl$ than $u^2$, at the cost of additional $O(n)$ term in the running time. 

This phenomenon compliments the better understood situation in which dynamical programming algorithms run faster in the {\em sparse} regime.  For instance, an upper bound on the support size of a solution is one of the parameters in the algorithms of Chan and He \cite{ChanHe} and Deng, Mao and Zhong \cite{DMZ}.

Let us choose a density level $\dl=\dl(n)$ with $2\le \dl\le n-1$ and consider the following scenario:

\begin{eqnarray}\begin{aligned}
&\mbox{An optimal solution to the problem }\eqref{initial_IP} \mbox{ is given by a vertex }{\ve z}^*\\&\mbox{ of the integer hull }\KP_I({\ve a},b)\mbox{ with }\|{\ve z}^*\|_0\ge \dl(n)+1\,.\label{initial_IP_dense}
\end{aligned}
\end{eqnarray}

One should note that \eqref{initial_IP_dense} imposes a lower bound on $\|{\ve a}\|_2$. By \eqref{t:Bound_via_determinant_bound1}, noting that $\Delta({\ve a}^\top)=\|{\ve a}\|_2$, we have
\begin{equation*} 
l \le \log_2\left(\frac{\|{\ve a}\|_2}{\gcd({\ve a})}\right)
\end{equation*}
and hence
\bea
\frac{\|{\ve a}\|_2}{\gcd({\ve a})}\ge2^{l-1}\,.
\eea

For convenience, we will use the notation 
\bea
L=\left\lfloor\frac{\dl(n)}{2^{\dl(n)}}\left\lfloor\frac{\|{\ve a}\|_2}{\gcd({\ve a})}\right\rfloor\right\rfloor\,,
\eea
that is $L$ is the integer part of the right-hand side in \eqref{cor_proximity_bound_order} with $d=\dl(n)$. 

Under assumption \eqref{initial_IP_dense}, we can reduce the position of the starting range $D[t-u+1, \ldots, t]$  to $t =uL$.

\begin{prop}
Suppose that $b>uL$. Then, assuming \eqref{initial_IP_dense}, we can compute $D[b]$ in $O(n)$ time from the array $D[uL-u+1, \ldots, uL]$.
\end{prop}
\begin{proof}
By Corollary  \ref{cor_EL_proximity_strong_ordered}, there exists an optimal solution 
${\ve z}^*$ with 
$
\|{\ve z}^*_{\tau}\|_1\le L\,.
$
Hence, 
\bea
{\ve a}_{\tau}{\ve z}^*_{\tau}\le uL\,.
\eea
Further, we have
\bea
{\ve a}_{\tau}{\ve z}^*_{\tau}+{\ve a}_{\bar\tau}{\ve z}^*_{\bar\tau}=b\,.
\eea
Since $0<{\ve a}_{\bar\tau}\le u$,  we have for some nonnegative integer $x< {\ve z}^*_{\bar\tau}$
\bea
y:={\ve a}_{\tau}{\ve z}^*_{\tau}+{\ve a}_{\bar\tau}x\in [uL-u+1, \ldots, uL]\,.
\eea
Without loss of generality,  we may assume that $y$ is the largest integer with this property.

We compute $n$ integers
\bea
y_i=b-a_i z_i\in [uL-u+1, uL]\,, \;z_i=\left\lceil \frac{b-uL}{a_i}\right\rceil
\eea
for $i\in [n]$ and notice that for at least one index $j$ we have  $y_j=y$.

Finally,
\bea
D[b]=\max_{i\in[n]} (D[y_i]+ c_iz_{i})\,.
\eea

\end{proof}

Note that 
$
\|{\ve a}\|_2\le \sqrt{n} \|{\ve a}\|_\infty= u\sqrt{n}\,,
$
which implies
\bea
t=uL\le 
\frac{u^2\dl(n)\sqrt{n}}{2^{\dl(n)}}\,.
\eea
Hence, we obtain an exponential  improvement in $\dl(n)$ on the classical bound $u^2$ for sufficiently large 
density level $\dl(n)$. In fact, an improvement already appears for the relatively small density levels comparable with $\log_2(n)$.

Indeed, if $\dl(n)>\lceil k\log_2(n)\rceil $ with some constant $k>1/2$, then
\bea
\frac{t}{u^2}=\frac{\dl(n)\sqrt{n}}{2^{\dl(n)}}\le \frac{\lceil k\log_2(n)\rceil }{n^{k-1/2}}=o(1)\,.
\eea

\section{Volumes and linear transforms}

In this section, we state a geometric lemma from \cite{ACH_IPCO,ACH_MPB}  needed for the proof of Theorem \ref{thm_EL_interger_hull}.

We begin by introducing some geometric notation. For a matrix $\UU\in \R^{l\times r}$, let $\lin(\UU)$ denote  the linear subspace of $\R^l$ spanned by the columns of $\UU$. Given a set $M\subset \R^r$, we define $\UU M=\{\UU{\ve x}\in \R^l: {\ve x}\in M\}$.   Thus, $\UU M$ is the image of $M$ under the linear map determined by $\UU$.
Next, let  $S$ be a $(l-k)$-dimensional subspace of $\R^l$. Choose an orthonormal basis ${\ve s}_1, \ldots, {\ve s}_{l-k}, {\ve s}_{l-k+1}, \ldots, {\ve s}_l$ of $\R^l$ such that the first $l-k$ vectors form a basis of $S$. Define $\SS_{l-k}=({\ve s}_1, \ldots, {\ve s}_{l-k})\in \R^{l\times (l-k)}$ and $\SS_{k}=({\ve s}_{l-k+1}, \ldots, {\ve s}_{l})\in \R^{l\times k}$. 

Let $\vol_i(\cdot)$  denote the $i$-dimensional volume. Given a Lebesgue measurable set  $M$ in the subspace $S$, we are interested in the $(l-k)$-dimensional volume of its image  $\DD M$, where $\DD\in \R^{l\times l}$ is assumed to be invertible. The following lemma provides a general expression for $\vol_{l-k}(\DD M)$ in terms of $\vol_{l-k}(M)$.

\begin{lemma}[Lemma 1 in \cite{ACH_IPCO,ACH_MPB}] \label{section_volume_regular} Let $S$ be an $(l-k)$-dimensional subspace of $\R^l$. Let $M$ be a Lebesgue measurable set in $S$, $\DD\in \R^{l\times l}$ be nonsingular
and let the rows of $\BB\in \R^{k\times l}$ form a basis of the subspace $(\DD S)^\bot$, the orthogonal complement of the subspace $\DD S=\lin(\DD \SS_{l-k})$. Then
\bea
\vol_{l-k}(\DD M)= |\det(\DD)| \sqrt{\frac{\det(\BB \BB^\top)}{\det(\BB\DD\DD^\top \BB^\top)}}\vol_{l-k}(M)\,.
\eea
\end{lemma}

\section{Proof of Theorem \ref{thm:UpperBound}}
First, we recall the definitions and lemmas used in the proof of  \cite[Theorem 3]{ACH_IPCO,ACH_MPB}.
Let ${\ve y} \in\R^n$ and define the open cube
\bea
C^n({\ve y})=\{{\ve x}\in  \R^n: \|{\ve x}-{\ve y} \|_{\infty} < 1 \}\,.
\eea
Given two points ${\ve u}, {\ve v}\in \R^n$, let
\bea
D({\ve u}, {\ve v})=\conv(C^n({\ve u}), C^n({\ve v}))\,.
\eea
That is, $D({\ve u}, {\ve v})$ is the convex hull of the two open cubes centered at ${\ve u}$ and ${\ve v}$.
\begin{lemma}[Lemma 3 in \cite{ACH_IPCO,ACH_MPB}]
\label{D_property} Let ${\ve u}, {\ve v}\in\R^n_{\ge 0}$. Then
$D({\ve u}, {\ve v})\cap\Z^n=D({\ve u}, {\ve v})\cap\Z^n_{\ge 0}\,.$
\end{lemma}
Next, let $E=E({\ve u}, {\ve v})$ denote the origin-symmetric open convex set defined by
\bea
E=\conv(C^n({\ve u}-{\ve v}), C^n({\ve v}-{\ve u}))\,.
\eea
Notice that
\be\label{zonotope_union}
E= (D({\ve u}, {\ve v})-{\ve v})\cup(-D({\ve u}, {\ve v})+{\ve v})\,.
\ee

\begin{lemma}[Lemma 4 in \cite{ACH_IPCO,ACH_MPB}]
\label{E_property} Suppose that ${\ve u},{\ve v}\in \KP(\AA, {\ve b})$. Then the bound 
\bea\label{E_property_bound} 
\vol_{n-m}(E\cap\ker(\AA))\ge 2^{n-m}(1+\|{\ve u}-{\ve v}\|_2)
\eea
holds.
\end{lemma}

We are now ready to prove Theorem \ref{thm:UpperBound}.
Choose an optimal vertex solution ${\ve x}^*$ of \eqref{initial_LP}. Clearly, we may assume  that ${\ve x}^*$ is not integer. We may also assume, without loss of generality, that there exists a vertex ${\ve z}^*$ of $\KP_I(\AA, {\ve b})$ that is the unique optimal solution to \eqref{initial_IP}. We justify this assumption below.

If ${\ve c}$ is orthogonal to $\ker(\AA)$, then we replace ${\ve c}$ with any vector $\bar{\ve c}$
such that ${\ve x}^*$ is an optimal vertex solution to 
\begin{eqnarray}
\max\{\bar\cc\cdot\xx:\xx\in\KP(\AA,\bb)\}\label{initial_LP_bar_c}
\end{eqnarray}
and the problem  
\begin{eqnarray}
\max\{\bar\cc\cdot\xx:\xx\in\KP(\AA,\bb)\cap\Z^{n}\}\label{initial_IP_bar_c}
\end{eqnarray}
has a unique optimal solution ${\ve z}^*$. Suppose now that ${\ve c}$ is not orthogonal to $\ker(\AA)$.  
In this case, without loss of generality, we may assume that ${\ve c}\in \ker(\AA)$. 

Let $V$ be the set of optimal solutions for \eqref{initial_IP}  that are vertices of $\KP_I(\AA, {\ve b})$.  
 By perturbing the original vector ${\ve c}$,  we may choose a vector $\bar{\ve c}$ such that ${\ve x}^*$ is the unique optimal solution to \eqref{initial_LP_bar_c} and the set $\overline{V}$ of optimal solutions to \eqref{initial_IP_bar_c} which are the vertices $\KP_I(\AA, {\ve b})$ is a subset of $V$. 
 
 If $\overline{V}$ consists of a single point, the claim is justified. Otherwise, for any point ${\ve z}^*\in \overline{V}$ the vector $\bar{\ve c}$ is on the boundary of the normal cone $N$ of $\KP_I(\AA, {\ve b})$ at the vertex ${\ve z}^*$, intersected with $\ker(\AA)$. Observe that, by uniqueness of ${\ve x}^*$ as the optimal solution to\eqref{initial_LP_bar_c},  the vector $\bar{\ve c}$ belongs to the relative interior of the normal cone  $N'$ of $\KP(\AA, {\ve b})$ at ${\ve x}^*$, intersected with the subspace $\ker(\AA)$. Therefore, perturbing  $\bar{\ve c}$, we may assume that it belongs to the relative interiors of both cones, $N\cap \ker(\AA)$ and $N'\cap \ker(\AA)$. Consequently, the vertex ${\ve z}^*$ of $\KP_I(\AA, {\ve b})$ is the unique optimal solution to \eqref{initial_IP_bar_c}, justifying the claim.

Suppose, to derive a contradiction, that the bound \eqref{t:Bound_via_pencils} does not hold for the point ${\ve z}^*$, that is
\be\label{lower_bound_distance}
\|{\ve x}^*-{\ve z}^*\|_2> \frac{\Delta(\AA)}{\gcd(\AA)}-1 = \det(\Lambda(\AA))-1\,.
\ee

The lower bound \eqref{lower_bound_distance} and Lemma \ref{E_property}  imply that for $E=E({\ve x}^*, {\ve z}^*)$
we have
\be\label{vol_E_kerA_bound}
\vol_{n-m}(E\cap\ker(\AA))> 2^{n-m} \det(\Lambda(\AA))\,.
\ee
Noting the bound \eqref{vol_E_kerA_bound}, Minkowski's first fundamental theorem (in the form of Theorem II in Chapter III of \cite{cassels1996introduction}) implies that the set $E\cap\ker(\AA)$
 contains nonzero points $\pm {\ve h}$ of the lattice $\Lambda(\AA)$. Using \eqref{zonotope_union}, we may assume without loss of generality
 that  ${\ve h}\in D({\ve x}^*, {\ve z}^*)-{\ve z}^*$. Consequently, the integer point ${\ve w}={\ve h}+{\ve z}^*$ belongs to the set $D({\ve x}^*, {\ve z}^*)\cap (\ker(\AA)+{\ve z}^*)$. By Lemma \ref{D_property}, we have ${\ve w}\in \Z^n_{\ge 0}$, and therefore ${\ve w}\in \KP(\AA, {\ve b})$.

We will show that 
 \be\label{objective_drops}
 \bar{\ve c}\cdot{\ve w}> \bar{\ve c}\cdot{\ve z}^*\,,
 \ee
contradicting the optimality of ${\ve z}^*$. 
 
Recall that ${\ve h}={\ve w}-{\ve z}^*$. Suppose that \eqref{objective_drops} does not hold. Since  ${\ve z}^*$ is the unique optimal solution to \eqref{initial_IP_bar_c}, it follows that 
\be\label{cr_negative}
\bar{\ve c}\cdot{\ve h}< 0\,.
\ee
We next show that 
for sufficiently small $\delta>0$, the point ${\ve x}^*-\delta{\ve h}$ is in the polyhedron $\KP(\AA, {\ve b})$.
Indeed, the point ${\ve x}^*-\delta{\ve h}$ is in the affine subspace defined by the system $\AA{\ve x}={\ve b}$. Thus, if ${\ve x}^*-\delta{\ve h}$ is not in $\KP(\AA, {\ve b})$ for every $\delta>0$, then there exists an index $j$ such that $x^*_j=0$ and $th_j>0$ for all $t>0$. 
Now, since $x^*_j=0$, every integer point ${\ve y}\in D({\ve x}^*, {\ve z}^*)$ satisfies $y_j\le z^*_j$ by construction of the set $D({\ve x}^*, {\ve z}^*)$. Further, $th_j>0$ for a positive $t$ implies $h_j>0$.
But then 
$w_j=h_j+z^*_j>z^*_j$, contradicting the inclusion ${\ve w}\in D({\ve x}^*, {\ve z}^*)$.

Finally, \eqref{cr_negative} implies that $\bar{\ve c}\cdot({\ve x}^*-\delta{\ve h})>\bar{\ve c}\cdot{\ve x}^*$, contradicting the optimality of 
the vertex ${\ve x}^*$.

\qed

\section{Proof of Theorem \ref{thm:UpperBoundOptimality}}

Let $k$ be a positive integer.  Define an $m\times n$ matrix $\AA(k)$ with columns $\AA_{\{1\}}(k)=-{\ve e}_1$ and $\AA_{\{n-m+1\}}(k)=k{\ve e}_1$, amended by the columns $\AA_{\{2\}}(k)=\cdots =\AA_{\{n-m\}}(k)={\ve e}_1$ and $\AA_{\{n-m+i\}}(k)={\ve e}_i$ for $i=2,\ldots, m$.  We also set ${\ve b}=(1,\ldots,1)^\top$.
Thus,
\bea
\AA(k)= \left(\begin{array}{rccccccc}
-1           & 1           & \cdots & 1         & k           & 0                & \cdots & 0  \\
 0           & 0           & \cdots & 0         &0             & 1                & \cdots & 0 \\
 \vdots   & \vdots    & \ddots & \vdots &\vdots     & \vdots         & \ddots & \vdots\\
0            & 0           & \cdots  & 0         &0            & 0                & \cdots  & 1\\
\end{array}
\right)\mbox{ and }
{\ve b}=\left(\begin{array}{c}
1\\
1\\
\vdots\\
1
\end{array}
\right)\,.
\eea
First, observe that $\KP(\AA(k), {\ve b})=\{{\ve x}\in\R^n: \AA(k){\ve x}={\ve b}\}\cap\R^n_{\ge 0}$. Hence, every nonempty intersection of $\KP(\AA(k), {\ve b})$ with a coordinate subspace of $\R^n$ is a face of the polyhedron $\KP(\AA(k), {\ve b})$. Moreover, by construction, every point of  $\KP(\AA(k), {\ve b})$ satisfies $x_{n-m+2}=\cdots=x_n=1$. 
 Consider the one-dimensional face $F$, that is an edge  of  $\KP(\AA(k), {\ve b})$, obtained by intersecting the polyhedron with the two-dimensional coordinate subspace spanned by ${\ve e}_{1}$ and ${\ve e}_{n-m+1}$, translated by  $x_{n-m+2}=\cdots=x_n=1$. After disregarding the zero variables, and the fixed variables $x_{n-m+2}, \ldots, x_n$, the face $F$ takes the form $F'=\{(x_1, x_{n-m+1})^\top\in \R^2_{\ge 0}: -x_1+kx_{n-m+1}=1\}$.  Observe that the point $(0,1/k)^\top$ is the unique vertex of $F'$, while the point $(k-1,1)^\top$ is the unique vertex of the integer hull $F'_I$. 
Therefore, the point ${\ve x}^*(k)=(0, \ldots, 0, 1/k, 1,\ldots, 1)^\top=(1/k){\ve e}_{n-m+1}+{\ve e}_{n-m+2}+\cdots+{\ve e}_{n}$ is the vertex of $F$, and, consequently, a vertex of $\KP(\AA(k), {\ve b})$.  Similarly, the point ${\ve z}^*(k)=(k-1, 0\ldots, 0, 1, 1,\ldots, 1)^\top=(k-1){\ve e}_{1}+{\ve e}_{n-m+1}+\cdots+{\ve e}_{n}$ is the vertex of $F_I$, and therefore a vertex of $\KP_I(\AA(k), {\ve b})$.  

Consequently, there exists an objective function vector ${\ve c}$ such that the original problem \eqref{initial_IP} with $\AA=\AA(k)$ and its linear relaxation \eqref{initial_LP} have unique optimal solutions ${\ve z}^*(k)$ and ${\ve x}^*(k)$, respectively. 

Using \eqref{CB}, we obtain
\bea
\Delta(\AA(k))=\sqrt{k^2+n-m}
\eea
and 
\bea
\|{\ve x}^*(k)-{\ve z}^*(k)\|_2= \sqrt{(k-1)^2+((k-1)/k)^2}\,.
\eea
Therefore,
\bea
\frac{\|{\ve x}^*(k)-{\ve z}^*(k)\|_2}{\Delta(\AA(k))-1}\rightarrow 1\mbox{ as } k\rightarrow \infty\,,
\eea
showing that \eqref{t:Bound_via_pencils_optimality} holds for sufficiently large $k$.

\begin{remark}
\label{worse_than_b}
As noted in the introduction, already for $m=1$ the minimum value of the product function at an optimal solution can be arbitrarily larger than $b$. Indeed, the optimal solution ${\ve z}^*(k)$ is unique, and for $m=1$ we have $\PF({\ve z}^*(k))=2k$, and $b=1$. 
\end{remark}

\section{Proof of Theorem \ref{thm_EL_interger_hull}}

Take any vertex ${\ve z}^*$ of $\KP_I(\AA, {\ve b})$.
For convenience, define $u_i=z_i^*+1$ for $i\in[n]$. Then $\DD=\DD({\ve z}^*)$ can be written as
$ \DD=\diag(u_1,\ldots, u_n)$ and the bound \eqref{EL_integer_hull_proximity_bound} becomes
\be\label{EL_product_bound}
u_1\cdots u_n \le \frac{\Delta(\AA \DD)}{\gcd(\AA)}=\ub_{\AA}({\ve z}^*)\,.
\ee
Let $\ker(\AA)$ denote the kernel subspace of the matrix $\AA$, that is $\ker(\AA)=\{{\ve x}\in\R^n: \AA{\ve x}={\ve 0}\}$.
Consider the box section
\bea
K=\ker(\AA)\cap (-u_1, u_1)\times\cdots\times (-u_n, u_n)\,.
\eea
There exists a subspace $S$ of $\R^n$ with $\dim(S)=\dim(\ker(\AA))= n-m$ such that
\bea
K=\DD Q\,,
\eea
where $Q=S\cap (-1,1)^n$ is a section of the open cube $(-1,1)^n$. 
Applying Lemma \ref{section_volume_regular}, we obtain
\be\label{EL_box_section_volume}
\vol_{n-m}(K)=
u_1\cdots u_n \frac{\Delta(\AA)}{\Delta(\AA\DD) } \vol_{n-m}(Q)\,.
\ee
By Vaaler's cube slicing inequality \cite{Vaaler1979}, 
\be\label{EL_Vaalers_bound_strong}
\vol_{n-m}(Q)\ge 2^{n-m}\,.
\ee 

Suppose, to derive a contradiction,  that \eqref{EL_product_bound} does not hold, that is,
\be\label{EL_contrary_product_bound}
u_1\cdots u_n >\frac{\Delta(\AA\DD)}{\gcd(\AA)}\,.
\ee
Let $\Lambda(\AA)=\ker(\AA)\cap\Z^n$ be the  lattice formed by all integer points in the kernel subspace of the matrix $\AA$. Combining \eqref{EL_box_section_volume}, \eqref{EL_Vaalers_bound_strong} and \eqref{EL_contrary_product_bound},
we obtain
\be\label{vol_K_lower_bound}
\begin{aligned}
\vol_{n-m}(K) > 2^{n-m}\frac{\Delta(\AA)}{\gcd(\AA)}=2^{n-m}\det(\Lambda(\AA))\,.
\end{aligned}
\ee
The $(n-m)$-dimensional subspace $\ker(\AA)$ can be considered as a usual Euclidean $(n-m)$-dimensional space. 
Hence, we may apply Minkowski's first fundamental theorem (in the form of Theorem II in Chapter III of \cite{cassels1996introduction}) to the set $K$ and the lattice $\Lambda(\AA)$.  Noting the bound \eqref{vol_K_lower_bound}, the  fundamental theorem implies that $K$ contains nonzero lattice  points $\pm {\ve h}\in \Lambda(\AA)$.

Consider the integer points ${\ve u}^+={\ve z}^*+ {\ve h}$ and ${\ve u}^-={\ve z}^*-  {\ve h}$.
By construction of $K$, both ${\ve u}^+$ and ${\ve u}^-$ have nonnegative coordinates. Hence, 
${\ve u}^+,{\ve u}^-\in \KP_I(\AA, {\ve b})$. 

Finally, ${\ve z}^* =({\ve u}^+ + {\ve u}^-)/2$, contradicting the assumption that ${\ve z}^*$ is a vertex of the integer hull $\KP_I(\AA, {\ve b})$.

\section{Proof of Theorem \ref{thm_EL_integer_hull_proximity_bound_optimality}}

We will reuse the construction from the proof of Theorem \ref{thm:UpperBoundOptimality}.
Recall that the integer hull $\KP_I(\AA(k), {\ve b})$ has a vertex ${\ve z}^*(k)=(k-1, 0\ldots, 0, 1, 1,\ldots, 1)^\top=(k-1){\ve e}_{1}+{\ve e}_{n-m+1}+\cdots+{\ve e}_{n}$. Hence,
\bea
\PF({\ve z}^*(k))=k 2^m\,.
\eea
Next, using \eqref{ub_expanded},
\bea
\begin{aligned}
\Delta(\AA(k)\DD({\ve z}^*(k))) &=\sqrt{\sum_{\substack{\sigma\subset[n]: \\ |\sigma|=m}}(\det(\AA_\sigma(k)) \prod_{i\in \sigma}(z^*_i(k)+1) )^2}\\
&=\sqrt{(k2^{m-1})^2+(n-m-1)(2^{m-1})^2+(k2^m)^2}\,.
\end{aligned}
\eea
It therefore suffices to observe that 
\bea
\frac{\PF({\ve z}^*(k))}{\Delta(\AA(k)\DD({\ve z}^*(k)))}\rightarrow\frac{2}{\sqrt{5}} \mbox{ as } k\rightarrow \infty\,.
\eea
Consequently, \eqref{EL_integer_hull_proximity_optimality} holds for sufficiently large $k$.

\section{Proof of Theorem \ref{thm_nonnegative}}

First we settle the case  $\|{\ve z}^*\|_0=n$.
Consider the zonotope
\bea
Z_{\AA}({\ve z}^*)=\sum_{i=1}^n[{\ve 0}, z^*_i \AA_{\{i\}}]\,.
\eea
Using \eqref{ub_expanded} and noting that all entries of ${\ve z}^*$ are positive, we obtain
\be\label{Delta_zono}
\begin{aligned}
\Delta(\AA\DD({\ve z}^*)) &=\sqrt{\sum_{\substack{\sigma\subset[n]: \\ |\sigma|=m}}(\det(\AA_\sigma) \prod_{i\in \sigma}(z^*_i+1) )^2}\\
&\le 2^m\sqrt{\sum_{\substack{\sigma\subset[n]: \\ |\sigma|=m}}(\det(\AA_\sigma) \prod_{i\in \sigma}z^*_i )^2}\\
&\le 2^m \sum_{\substack{\sigma\subset[n]: \\ |\sigma|=m}}|\det(\AA_\sigma) \prod_{i\in \sigma}z^*_i|=2^m\vol_m(Z_{\AA}({\ve z}^*))
\,,
\end{aligned}
\ee
where for the last equality we refer the reader to Shephard \cite[Equation (57)]{Shephard} (see also Gover and Krikorian \cite[Corollary 3.4]{Gover_Krikorian}). 

Since ${\ve b}=\sum_{i=1}^{n}z^*_i \AA_{\{i\}}$, the entries $z^*_i$ are positive, the columns $\AA_{\{i\}}$ are nonnegative, and $\AA$ has rank $m$, it follows that  ${\ve b}\in\Z^m_{>0}$ and 
\bea
Z_{\AA}({\ve z}^*)\subset [0,b_1]\times \cdots \times [0,b_m]\,.
\eea
Consequently,
\be\label{zono_vol}
\vol_m(Z_{\AA}({\ve z}^*))\le \prod_{i=1}^m b_i\,.
\ee
Combining \eqref{Delta_zono}, \eqref{zono_vol} and \eqref{EL_integer_hull_proximity_bound} yields \eqref{nonnegative_A}.

We now consider the case $s=\|{\ve z}^*\|_0<n$. 
Let $\tau=\supp({\ve z}^*)$, and consider the subvector ${\ve y}^*={\ve z}^*_\tau$ together with the submatrix $\AA_\tau$.

If $\AA_\tau$ has rank $r<m$, then we disregard $m-r$ rows with indices in a set $\mu=\{i_1,\ldots, i_{m-r}\}$ so as to obtain an $r\times s$ submatrix $\widetilde\AA$ of $\AA_\tau$ of rank $r$. 
Let also $\tilde{\ve b}$ denote the subvector obtained from ${\ve b}$ by deleting the entries with indices $i_1,\ldots, i_{m-r}$.
If $\AA_\tau$ has rank $r=m$, we simply set $\widetilde\AA=\AA_\tau$ and $\tilde{\ve b}={\ve b}$.

The polyhedron $\KP(\widetilde\AA, \tilde{\ve b})$ is a face of $\KP(\AA, {\ve b})$. Hence ${\ve y}^*$ is a vertex of the integer hull $\KP_I(\widetilde\AA, \tilde{\ve b})$. 

If $\|{\ve y}^*\|_0=s>r$, we may apply the above argument with $\AA=\widetilde\AA$,  ${\ve b}=\tilde{\ve b}$ and 
${\ve z}^*={\ve y}^*$. Consequently,
\bea
\PF({\ve y}^*)\le 2^r  \prod_{i\notin \mu} \tilde b_i\le 2^m \prod_{i\in\supp({\ve b})}  b_i\,.
\eea
It remains to consider the case $s=r$. In this case, we consider the parallelepiped 
\bea
Q=\sum_{i\in\tau}[{\ve 0}, y^*_i \widetilde\AA_{\{i\}}]
\eea
and observe that $Q$ is contained in the box $[0,\tilde b_{k_1}]\times \cdots \times [0, \tilde b_{k_r}]$,  where $\{k_i,\ldots, k_r\}=[m]\setminus \mu$. Therefore,
\bea
\PF({\ve y}^*) |\det(\widetilde\AA)| \le 2^s |\det(\widetilde\AA)| \prod_{i\in \tau} {y}_i^* = 2^s \vol_r(Q) \le  2^m \prod_{i\in\supp({\ve b})}  b_i\,,
\eea
which completes the proof.

\section{Proof of Theorem \ref{thm_EL_proximity_strong_ordered}}

Choose a vertex ${\ve z}^*$ of the integer hull $\KP_I(\AA, {\ve b})$. By Theorem \ref{thm_EL_interger_hull}, the point ${\ve z}^*$ satisfies \eqref{EL_integer_hull_proximity_bound}.
 %
Applying \eqref{ub_expanded}, we obtain
\bea
\Delta(\AA\DD({\ve z}^*))=\sqrt{\sum_{\substack{\sigma\subset[n]: \\ |\sigma|=m} }(\det(\AA_\sigma) \prod_{i\in \sigma}(z_i^*+1) )^2}\,.
\eea

Hence, by \eqref{EL_integer_hull_proximity_bound}, 
\be\label{EL_integer_hull_proximity_bound_working_copy_1} 
\PF({\ve z}^*)=\prod_{i\in [n]} (z_i^*+1)\le (\gcd(\AA))^{-1}\sqrt{\sum_{\substack{\sigma\subset[n]: \\ |\sigma|=m} }(\det(\AA_\sigma) \prod_{i\in \sigma}(z_i^*+1) )^2}\,.
\ee

Recall that  $z^*_{i_1}\le z^*_{i_2}\le \cdots \le z^*_{i_n}$ and $\tau=\{i_{1}, \ldots, i_{n-m}\}$. Divide both sides of  
the inequality \eqref{EL_integer_hull_proximity_bound_working_copy_1} by the product $(z^*_{i_{n-m+1}}+1)\cdots (z^*_{i_n}+1)$, constructed using  the $m$ largest entries of ${\ve z}^*$. We obtain the desired bound
\bea 
\begin{aligned}
\PF({\ve z}^*_\tau) =\prod_{i\in \tau} (z_i^*+1)&\le (\gcd(\AA))^{-1}\sqrt{\sum_{\substack{\sigma\subset[n]: \\ |\sigma|=m} }(\det(\AA_\sigma))^2}\\
&=\frac{\Delta(\AA)}{\gcd(\AA)}\,,
\end{aligned}
\eea
where the final equality follows from \eqref{CB}.

\section{Proof of Theorem \ref{thm_EL_proximity_strong_ordered_optimality}}

Let $k> 1$ be an integer and set $n=m+k-1$. We consider the $m\times n$ matrix
\bea
\AA(k)= \left(\begin{array}{llccccc}
2^{k-1}           & 2^{k-2}           & \cdots & 1          & 0                & \cdots & 0  \\
 0           & 0           & \cdots & 0                     & 1                & \cdots & 0 \\
 \vdots   & \vdots    & \ddots & \vdots      & \vdots         & \ddots & \vdots\\
0            & 0           & \cdots  & 0                   & 0                & \cdots  & 1\\
\end{array}
\right)\mbox{ and }
{\ve b}(k)=\left(\begin{array}{c}
2^k-1\\
1\\
\vdots\\
1
\end{array}
\right)\,.
\eea

For convenience, we denote by ${\ve 1}_{n}=(1,\ldots,1)^\top$ the  $n$-dimensional all-ones vector. We will begin with the following lemma. 
\begin{lemma}
The point ${\ve 1}_{n}$ is a vertex of the integer hull $\KP_I(\AA(k), {\ve b}(k))$.
\end{lemma}
\begin{proof}
We proceed by induction on $k$. For the basis step $k=2$, the claim holds since the only  integer points in the polytope $\KP(\AA(2), {\ve b}(2))$ are ${\ve 1}_n$ and $(0,3,1,\ldots, 1)^\top$. 

Assume now that the claim does not hold for some $k\ge 3$. Observe that any integer point ${\ve z}=(z_1, \ldots, z_n)^{\top}\in \KP(\AA(k), {\ve b}(k))$ satisfies $z_1\le 1$. Hence ${\ve 1}_n$ belongs to the face $F=\KP_I(\AA(k), {\ve b}(k))\cap\{{\ve x}\in \R^n: x_1=1\}$ of the integer hull $\KP_I(\AA(k), {\ve b}(k))$. Since ${\ve 1}_n$ is not a vertex, it must be a convex combination of integer points in $F$. Next,  observe that, in the variables $x_2, \ldots, x_n$, the face $F$ coincides with $\KP_I(\AA(k-1), {\ve b}(k-1))$. This implies that  ${\ve 1}_{n-1}$ is not a vertex of $\KP_I(\AA(k-1), {\ve b}(k-1))$, contradicting the induction hypothesis. This completes the proof.
\end{proof}

Using \eqref{ub_expanded}, we obtain
\bea
\Delta(\AA(k))=\sqrt{2^{2(k-1)}+2^{2(k-2)}+\cdots+1}\,.
\eea
Taking ${\ve z}^*={\ve 1}_n$ and choosing any $\tau$ as in the statement of Theorem \ref{thm_EL_proximity_strong_ordered}, we get
\bea
\PF({\ve z}^*_\tau)= 2^{n-m}=2^{k-1}\,.
\eea
Finally, using the geometric series formula, we have
\bea
\frac{\PF({\ve z}^*_\tau)}{\Delta(\AA(k))}=\left(\sum_{i=0}^{k-1} 2^{-2i}\right)^{-1/2}=\left(\frac{4}{3}\left(1- \frac{1}{4^k}\right)\right)^{-1/2}\rightarrow \frac{\sqrt{3}}{2}\mbox{ as } k\rightarrow \infty\,
\eea
and therefore \eqref{HW_proximity_bound_order_optimality} holds for sufficiently large $k$.

\begin{remark}
\label{tightness_b_plus_1}
Notice that for $m=1$ we have $\PF({\ve z}^*)=2^k=b+1$. Hence the bound \eqref{DMZ-bound} is tight.
\end{remark}

\section{Proof of Corollary \ref{cor_general_sparsity}}

We can rewrite \eqref{HW_proximity_bound_order} as
\bea
2^{\|{\ve z}^*_{\tau}\|_0}\le \frac{2^{\|{\ve z}^*_{\tau}\|_0}\Delta(\AA)}{\PF({\ve z}^*_\tau)\gcd(\AA)} \,.
\eea
Taking the logarithm with base 2 of both sides, and using the observation that 
 $\|{\ve z}^*\|_0\le \|{\ve z}^*_{\tau}\|_0+m$, we obtain the bound \eqref{delicate_general_sparsity}.

\section{Proof of Corollary \ref{cor_EL_proximity_strong_ordered}}

We will use the following elementary lemma proved in Aliev et al. \cite{ACHW}.
\begin{lemma}[Lemma 9 in \cite{ACHW}]\label{t:sum_product}
For $d\ge 2$ and $x_1, \ldots, x_d\ge 1$ the inequality
\bea\label{sum_product}
x_1+\cdots+x_d \le \frac{d(x_1+1)\cdots(x_d+1)}{2^{d}}
\eea
holds. \end{lemma}
Applying Lemma \ref{t:sum_product} to $d=\|{\ve z}^*_{\tau}\|_0$ nonzero entries of  ${\ve z}^*_{\tau}$, and using the bound \eqref{HW_proximity_bound_order}, we immediately obtain \eqref{cor_proximity_bound_order}.

\section{Proof of Theorem \ref{thm_EL_proximity_strong}}

Let  $\ve x^*$ be a  vertex optimal solution to the linear relaxation \eqref{initial_LP} with a basis $\gamma$.
Let ${\ve z}^*$ be a vertex of the integer hull $\KP_I(\AA, {\ve b})$ as in Theorem \ref{thm:UpperBound}. Then the points $\ve x^*$ and ${\ve z}^*$ satisfy \eqref{t:Bound_via_pencils}. Moreover, ${\ve z}^*$ also satisfies \eqref{EL_integer_hull_proximity_bound} by Theorem \ref{thm_EL_interger_hull}. 

Proceeding as in the proof of Theorem \ref{thm_EL_proximity_strong_ordered}, we obtain the bound \eqref{EL_integer_hull_proximity_bound_working_copy_1}.
Dividing both sides of  
the inequality \eqref{EL_integer_hull_proximity_bound_working_copy_1} by the product $(z^*_{i_{n-m+1}}+1)\cdots (z^*_{i_n}+1)$, constructed using the $m$ largest entries of ${\ve z}^*$, yields
\be\label{EL_integer_hull_proximity_divided} 
\begin{aligned}
\frac{\prod_{i\in [n]} (z_i^*+1)}{(z^*_{i_{n-m+1}}+1)\cdots (z^*_{i_n}+1)} &\le \sqrt{\sum_{\substack{\sigma\subset[n]: \\ |\sigma|=m} }(\det(\AA_\sigma))^2}\\
&=\Delta(\AA)\,.
\end{aligned}
\ee

Next, denoting for convenience $\sigma=\bar\tau=\{i_{n-m+1},\ldots, i_n\}$, we can write
\be\label{EL_integer_hull_proximity_ratio} 
\begin{aligned}
\frac{\prod_{i\in [n]} (z^*_i+1)}{(z^*_{i_{n-m+1}}+1)\cdots (z^*_{i_n}+1)} &= \frac{\prod_{i\in \gamma} (z^*_i+1)\prod_{i\in\bar \gamma} (z^*_i+1)}{(z^*_{i_{n-m+1}}+1)\cdots (z^*_{i_n}+1)} \\
&=\frac{\prod_{i\in \gamma, i\notin\sigma} (z^*_i+1)\prod_{i\in\bar \gamma} (z^*_i+1)}{\prod_{i\in \sigma\cap\bar \gamma} (z^*_i+1)}\,.
\end{aligned}
\ee
Recall that $\gamma$ is a basis of the vertex ${\ve x}^*$, so that ${\ve x}^*_{\bar\gamma}={\ve 0}$. Let ${\ve r}={\ve x}^*-{\ve z}^*$. Then for all $i\in\sigma\cap\bar \gamma$, we have $z^*_i=|x^*_i-z^*_i|=|r_i|$.  Applying the AM--GM inequality gives
\bea
\prod_{i\in \sigma\cap\bar \gamma} (z^*_i+1)= \prod_{i\in \sigma\cap\bar \gamma}(|r_i|+1)\le \left(\frac{\|{\ve r}\|_1}{m}+1\right)^m\,.
\eea
Using the Cauchy-Schwarz inequality, we further obtain
\bea
\prod_{i\in \sigma\cap\bar \gamma} (z^*_i+1)\le \left(\frac{\|{\ve r}\|_2}{\sqrt{m}}+1\right)^m\le \left(\frac{\Delta(\AA)-1}{\sqrt{m}}+1\right)^m\,,\eea
where the last inequality follows from \eqref{t:Bound_via_pencils}.
Next, using \eqref{EL_integer_hull_proximity_ratio}, we have
\be\label{EL_integer_hull_proximity_ratio_1} 
\begin{aligned}
\frac{\prod_{i\in [n]} (z^*_i+1)}{(z^*_{i_{n-m+1}}+1)\cdots (z^*_{i_n}+1)} &= \frac{\prod_{i\in \gamma, i\notin\sigma} (z^*_i+1)\prod_{i\in\bar \gamma} (z^*_i+1)}{\prod_{i\in \sigma\cap\bar \gamma} (z^*_i+1)}\\
&\ge \prod_{i\in\bar \gamma} (z^*_i+1)\left(\frac{\Delta(\AA)-1}{\sqrt{m}}+1\right)^{-m}\,, \end{aligned}
\ee

Combining the inequalities \eqref{EL_integer_hull_proximity_divided} and \eqref{EL_integer_hull_proximity_ratio_1}, we obtain the desired bound
\bea\label{EL_proximity_bound} 
\PF({\ve z}_{\bar\gamma}^*)=\prod_{i\in \bar\gamma} (z^*_i+1)\le \Delta(\AA)\left(\frac{\Delta(\AA)-1}{\sqrt{m}}+1\right)^m\,.
\eea

\end{document}